\documentclass[12pt]{amsart}
\usepackage{amssymb}
\usepackage{epsfig}
\usepackage{graphicx}

\theoremstyle{plain}
\newtheorem{prop}{Proposition}[section]
\newtheorem{theo}[prop]{Theorem}
\newtheorem{coro}[prop]{Corollary}

\theoremstyle{definition}

\newtheorem{conj}[prop]{Conjecture}
\newtheorem{thes}[prop]{Thesis}

\newtheorem{rema}[prop]{Remark}

\newtheorem{exam}[prop]{Example}

\newcommand{\ch}{\mathrm{ch}}
\newcommand{\td}{\mathrm{td}}

\newcommand{\Hom}{\mathrm{Hom }}
\newcommand{\Ext}{\mathrm{Ext }}

\newcommand{\Pic}{\mathrm{Pic}}
\newcommand{\Alb}{\mathrm{Alb}}

\newcommand{\NE}{\mathrm{NE}}
\newcommand{\oNE}{\overline{\mathrm{NE}}}

\newcommand{\rN}{\mathrm{N}}
\newcommand{\rH}{\mathrm{H}}

\newcommand{\Gr}{\mathrm{Gr}}

\newcommand{\ra}{\rightarrow}

\newcommand{\bP}{{\mathbb P}}

\newcommand{\bC}{{\mathbb C}}

\newcommand{\bL}{{\mathbb L}}
\newcommand{\bN}{{\mathbb N}}
\newcommand{\bQ}{{\mathbb Q}}
\newcommand{\bR}{{\mathbb R}}

\newcommand{\bZ}{{\mathbb Z}}

\newcommand{\cE}{{\mathcal E}}
\newcommand{\cF}{{\mathcal F}}
\newcommand{\cI}{{\mathcal I}}

\newcommand{\cO}{{\mathcal O}}

\newcommand{\fS}{{\mathfrak{S}}}
\newcommand{\sg}{{\mathsf{g}}}

\newcommand{\Def}{\mathrm{Def}}

\author{Brendan Hassett and Yuri Tschinkel}
\title[Intersection numbers and extremal rays]{Intersection numbers
of extremal rays on holomorphic symplectic varieties}

\begin{document}
\date{\today}

\maketitle

\section{Introduction} \label{sect:intro}

Suppose $X$ is a smooth projective complex variety.  
Let $\rN_1(X,\bZ) \subset \rH_2(X,\bZ)$ and $\rN^1(X,\bZ) \subset \rH^2(X,\bZ)$
denote the group of curve classes modulo homological equivalence and the 
N\'eron-Severi group respectively.  The monoids of effective classes
in each group generate cones $\NE_1(X) \subset \rN_1(X,\bR)$ and 
$\NE^1(X) \subset \rN^1(X,\bR)$ with closures $\oNE_1(X)$ and $\oNE^1(X)$,
the {\em pseudoeffective cones}.
These play a central r\^ole in the birational geometry of $X$.  
By Kleiman's criterion \cite{Kl}, a divisor $D$ is ample if and only if $D\cdot C>0$
for each $C\in \oNE_1(X)\setminus \{0\}$.  More generally, a divisor $D$ is
{\em nef} if $D\cdot C \ge 0$ for each such $C$.  Divisors in the interior of $\NE^1(X)$
are {\em big} and induce birational maps from $X$.  In general, very little
is known about the structure of these cones beyond Mori's Cone Theorem \cite[\S 3.3]{CKM} and 
more general finiteness results arising from the Minimal Model Program \cite[1.1.9]{BCHM}.

Let $X$ be an irreducible holomorphic symplectic variety, i.e.,
a smooth projective simply-connected manifold admitting a unique
nondegenerate holomorphic two-form.  
Let $\left(,\right)$ denote the Beauville-Bogomolov form on the 
cohomology group $\rH^2(X,\bZ)$, normalized so that it is integral
and primitive.  Duality gives a $\bQ$-valued form on $\rH_2(X,\bZ)$,
also denoted $\left(,\right)$. 
When $X$ is a K3 surface
these coincide with the intersection form.

For a polarized K3 surface $(S,g)$, we can 
read off the ample cone from the Hodge structure on the middle cohomology
and the polarizing class.  Precisely, the cone of effective curves
can be expressed \cite{LP}
$$\NE_1(S)=\left<C \in \rN_1(S,\bZ): \left(C,C\right) \ge -2, C\cdot g>0 \right>.$$
(We use the notation $\left<\cdots \right>$ for the `cone generated by the enclosed elements'.)  
Thus the ample cone can be expressed
$$\left<h \in \rN^1(S,\bZ):h\cdot C >0, \forall C\in \rN_1(S,\bZ) \text{ with }
	\left(C,C\right) \ge -2, C\cdot g >0 \right>.$$
Our goal is to extend these explicit descriptions to higher dimensions.  
Of course, on surfaces curves and divisors are equivalent.  In higher dimensions
we can seek characterizations of both the 
effective curves and the effective
divisors.  

\begin{thes} \label{thesis:main}
Let $X$ be an irreducible holomorphic symplectic manifold of dimension $2n$.
There is a universal constant $c_X \in \bQ_{>0}$ depending only on the
deformation equivalence class of $X$, with the following properties:
\begin{itemize}
\item
Let $X'$ be a deformation of $X$ containing a Lagrangian 
submanifold $P' \subset X'$ with $P'\simeq \bP^n$.
Let $\ell$ denote the class of a line on $P'$
interpreted as an element of $\rH_2(X,\bZ) \subset \rH^2(X,\bZ) \otimes_{\bZ} \bQ$.  
Then $\left(\ell,\ell\right)= -c_X$. 
\item
Assume that $X$ admits a polarization $g$.  Then we have
$$\NE_1(X) = \left< C \in \rN_1(X,\bZ): \left(C,C\right) \ge  -c_X, C\cdot g>0 \right>.$$
\end{itemize}
\end{thes}
Our main objective is to collect evidence supporting this thesis for various deformation equivalence
classes of irreducible holomorphic symplectic manifolds.  In doing so, we formulate
precise conjectures in these cases.  

It is natural to seek to relate the intersection properties of extremal rays $R \in \rN_1(X,\bZ)$
to the geometry of the associated contractions.  We carry out this analysis for
divisorial contractions in Section~\ref{sect:DC}, and obtain new structural 
results on cones of effective divisors.  

\subsection*{Conventions}
Let $R \in \rN_1(X,\bZ)$ denote the primitive generator of an extremal ray
on $\oNE_1(X)$ with negative self-intersection.  Let $\rho \in \rN^1(X,\bZ)$
denote the divisor class corresponding to $\rho$, i.e., using the form
$\left(,\right)$ to embed
$$\rH^2(X,\bZ) \subset \rH_2(X,\bZ)$$
we take $\rho$ to be the smallest positive multiple of $R$ contained
in $\rH^2(X,\bZ)$.  

The notation $S \sim M$ means `$S$ is deformation equivalent to $M$'.

\subsection*{Punctual Hilbert schemes}
Let $X$ be deformation equivalent to the punctual Hilbert scheme $S^{[n]}$,
where $S$ is a K3 surface. 
For $n>1$ the Beauville-Bogomolov form can be written
\cite[\S 8]{beauville}
$$\rH^2(X,\bZ) \simeq \rH^2(S,\bZ)_{\left(,\right)} \oplus_{\perp} \bZ \delta, \quad
				\left(\delta,\delta\right)=-2(n-1)$$
where $2\delta$ is the class of the `diagonal' divisor 
$\Delta^{[n]} \subset S^{[n]}$ parametrizing
non-reduced subschemes.  For each homology class $f\in \rH^2(S,\bZ)$, 
let $f \in \rH^2(X,\bZ)$ denote the class parameterizing subschemes with
some support along $f$.  This is compatible with the lattice embedding
above.  
Duality gives a $\bQ$-valued form 
on homology
$$
\rH_2(X,\bZ) \simeq \rH_2(S,\bZ)_{\left(,\right)} \oplus_{\perp} \bZ \delta^{\vee}, \quad
			\left(\delta^{\vee},\delta^{\vee}\right)=-1/2(n-1)$$
where $\delta^{\vee}$ is characterized as the homology class orthogonal
to $\rH^2(S,\bZ)$ and satisfying $\delta^{\vee}\cdot \delta =1$.  

Our thesis takes the following form:
\begin{conj} \label{conj:main1}
Let $(X,g)$ be a polarized variety deformation equivalent to $S^{[n]}$ 
where $S$ is a K3 surface.  
Then 
\begin{equation} \label{eqn:maincone}
\NE_1(X) = \left<R \in \rN_1(X,\bZ): \left(R,R\right) \ge -(n+3)/2, R \cdot g >0 \right>
\end{equation}
and a divisor class $h$ on $X$ is ample if
$h\cdot R > 0$ for each $R \in \rN_1(X,\bZ)$
satisfying $g\cdot R>0$ and 
$$\left(R,R \right) \ge -(n+3)/2.$$
\end{conj}
This generalizes to higher dimensions conjectures of \cite{HTGAFA99} in the
special case $n=2$;  see \cite{HTGAFA08} for a proof in this case 
that the divisors predicted to be
nef are in fact nef.  

For the examples in 
dimensions $\le 8$ discussed below, the geometry of each
extremal ray is characterized by its intersection properties:
Let $R$ be a primitive integral class with 
$\left(R,R\right) <0$
generating an extremal ray of the cone $\oNE_1(X)$. 
Choose $k\in \bN$ such that
$$ -(k+3)/2 \le  \left(R,R\right) <
		\begin{cases} -(k+2)/2 & \text{ if } k>1 \\
				0 & \text{ if }k=1.
			\end{cases}
$$
Then there exists an extremal contraction $\beta:X\ra Y$ of $R$ 
$$\begin{array}{ccccccc} 
\bP^k & \ra & Z_{\circ}  & \subset & Z & \subset & X \\
      &     &{\scriptstyle \beta_{\circ}} \downarrow \quad &         & \downarrow & & {\scriptstyle \beta} \downarrow  \quad  \\
      &     & W_{\circ} & \subset & W  &  \subset & Y 
\end{array}
$$
where $Y$ is the image, $Z$ the exceptional locus, $W$ its image,
and $W_{\circ} \subset W$ is an open subset over which $\beta_{\circ}$ is smooth.  
Furthermore, $\beta_{\circ}$ is a $\bP^k$-bundle over $W_{\circ}$.
Theorem~\ref{theo:divisorcone}, for general holomorphic symplectic manifolds, shows that divisorial
extremal rays satisfy $\left(R,R\right) \ge -2$.  

It is a general property of birational contractions of holomorphic symplectic
manifolds \cite[Theorem 2.3]{Kal} (cf. Namikawa \cite{Nam})
that $W_{\circ}$ is holomorphic symplectic and $W\setminus W_{\circ}$
has even codimension in $Y$.
Kaledin 
shows more:  He gives a finite `Poisson stratification'
of $W$.

\begin{rema}
For $n=5$, the geometric characterization sketched above is incompatible with Conjecture~\ref{conj:main1}.  
There are elements $R \in \rH_2(S^{[5]},\bZ)$ with 
$\left(R,R\right)=-9/8$;  however, a result of Markman implies these cannot correspond to divisorial contractions.
See Remark~\ref{rema:Markman} for details.  
\end{rema}

Finally, the results of \cite[\S 4]{HTGAFA99} (building on \cite{Ran}) imply that these conjectures are stable under small deformations:
Given a negative extremal ray $R=[\ell]$ for $X$ with an interpretation as above, for any small deformation
$X'$ of $X$ such that the class $R$ remains algebraic, there exists an $\ell' \subset X'$ specializing
to $\ell$ and subject to the same geometric interpretation.

\subsection*{Generalized Kummer varieties}
Suppose that $X$ is a holomorphic symplectic variety of dimension $2n$,
deformation equivalent to a generalized Kummer variety $K_n(A)$, defined as follows:  Given
an abelian surface $A$ with degree $(n+1)$ Hilbert scheme $A^{[n+1]}$, 
$K_n(A)$ is the fiber over $0$ of the addition map $\alpha:A^{[n+1]} \ra A$.
Here the Beauville-Bogomolov form is
$$\rH^2(K_n(A),\bZ) = \rH^2(A,\bZ) \oplus_{\perp} \bZ e, \quad \left(e,e\right)=-2(n+1)$$
where $2e$ is the class of the non-reduced subschemes
\cite[\S 4.3.1]{Yosh01}.  (Each class $\eta \in \rH^2(A,\bZ)$ yields a class in $\rH^2(K_n(A),\bZ)$, i.e.,
the subschemes with some support along $\eta$.)  As in the Hilbert scheme case, we use $\left(,\right)$ to
embed $\rH^2(X,\bZ) \subset \rH_2(X,\bZ)$ and extend $\left(,\right)$ to a $\bQ$-valued
form on $\rH_2(X,\bZ)$.  

We offer a precise statement analogous to the conjectures in \cite{HTGAFA99} for varieties deformation equivalent to 
$S^{[2]}$: 

\begin{conj} \label{conj:Kummer2}
Let $(X,g)$ be a polarized variety deformation equivalent to a four-dimensional
generalized Kummer variety.  
Then 
$$\NE_1(X)=\left<R \in \rN_1(X,\bZ): \left(R,R\right) \ge -3/2, R \cdot g >0 \right>$$
and a divisor class $h$ on $X$ is ample if and only if
$h\cdot R > 0$ for each $R \in \rN_1(X,\bZ)$
satisfying $g\cdot R>0$ and 
$$\left(R,R \right) \ge -3/2.$$
Furthermore, extremal rays $R$ with $\left(R,R\right)<0$ have the following interpretations:
\begin{itemize}
\item{$\left(R,R\right)=-\frac{1}{6}$:  There exists a divisor $E\subset X$ with class $[E]=2e$
where $e=6R$, with $E$ ruled over an abelian surface with fibers of class $R$.}
\item{$\left(R,R\right)=-\frac{2}{3}$:  There exists a divisor $F \subset X$ where $F=3R$,
with $F$ ruled over an abelian surface with fibers of class $R$.}
\item{$\left(R,R\right)=-\frac{3}{2}$:  There exists a plane $P \subset X$ whose lines $\ell \subset P$
have class $R$.}
\end{itemize}
\end{conj}
In particular, if $\rho$ is the smallest positive multiple of $R$ contained in $\rH^2(X,\bZ)$
(via the inclusion $\rH^2(X,\bZ) \subset \rH_2(X,\bZ)$) then in each case
$$\left(\rho,\rho\right)=\left(e,e\right)=-6.$$
In \cite{HT10}, we prove this last statement when $X$ contains a Lagrangian plane $P$
and $\rho$ is associated with the class of a line $\ell \subset P$;
if $R=[\ell]$ is primitive 
then $\left(R,R\right)=-\frac{3}{2}$.

{\bf Acknowledgments:} We are grateful to L. G\"ottsche for useful conversations. 
We are grateful to D. Kaledin for
reading an early draft of this manuscript, and to E. Markman for
pointing out the results of Boucksom and Druel we utilize.  The first
author was supported by National Science Foundation Grants 0554491 and 0901645;
the second author was supported by National Science Foundation
Grants 0554280 and 0602333.

\section{Divisorial contractions} \label{sect:DC}
Let $X$ be an irreducible holomorphic symplectic manifold of dimension $2n$.
What is needed to prove
that the divisors predicted to be ample are indeed ample?  This would entail showing that
all negative extremal rays $R$ on $X$ satisfy
$\left(R,R\right) \ge -c_X$.  In this section, we establish some uniform bounds of this
type, and offer applications to cones of moving and effective divisors.  

Let $\beta:X \ra Y$ be an extremal 
contraction of the ray $R$.  Let $Y_{sing}$ denote the singular locus of $Y$;
it is a result of Namikawa \cite[1.4]{Nam} that $Y_{sing}$ is even dimensional.
Moreover, Kaledin \cite[2.3]{Kal} has shown
that $Y_{sing}$ is stratified by symplectic manifolds, and $\beta$ has fiber dimension
$n-j$ over the $(2j)$-dimensional stratum.

Thus $\beta$ is divisorial if and only if $Y_{sing}$ has codimension two.
In this situation
\cite[\S 1]{Nam} and \cite[5.1]{Wier} offer detailed structural information about the singularities
at the generic point of $Y_{sing}$, i.e., they are rational double points of types
$A_1$ or $A_2$.  Let $E$ be the exceptional divisor and $C$ the class of an exceptional
rational curve.
Assume that $C$ is the class of an irreducible component $\ell$ of the fiber
of $\beta$ over a generic point of $Y_{sing}$.
First, we have
$$E\cdot C = \begin{cases} -2& \text{ in the $A_1$ case} \\
			   -1& \text{ in the $A_2$ case.}
		\end{cases}
$$
Regard the classes of $E$ and $C \in \rH_2(X,\bQ)=\rH^2(X,\bQ)$; the class of $E$ is a positive rational multiple
of $C$.  Indeed, we may deform $X$ preserving the algebraicity of the class of $C$;  the curve $C$ deforms to a curve
in nearby fibers \cite[\S 4]{HTGAFA99}, which still sweeps out a divisor that specializes back to $E$.

Express $E=m\rho$ and $C=m'R$  for $\rho\in \rH^2(X,\bZ)$ and $R\in \rH_2(X,\bZ)$ primitive and $m,m' \in \bN$.
We have $\rho=kR$ for some $k\in \bN$ as well, since the Beauville-Bogomolov form induces an inclusion
$\rH^2(X,\bZ) \subset \rH_2(X,\bZ)$.  Thus we have
$$E\cdot C = mm' \rho\cdot R=mm'k \left(R,R\right),$$
which implies that $\left(R,R\right) \ge -2$.  

We summarize this discussion in the following theorem:
\begin{theo} \label{theo:divisorial}
Let $X$ be an irreducible holomorphic symplectic manifold.  Let $R$ be a primitive generator of an
extremal ray of $X$ such that the associated extremal contraction is divisorial.  Then
$$-2 \le \left(R,R\right)  < 0.$$
\end{theo}

\begin{rema} \label{rema:footnote}
Our argument works under slightly weaker assumptions:  It suffices that
$\beta:X \ra Y$ be a divisorial contraction with irreducible exceptional
divisor $E$, and that $R$ be the primitive class associated with an irreducible
component of the generic exceptional fiber.  
Since $E$ is irreducible, the monodromy action on the irreducible components
of the generic exceptional fiber is transitive.  This suffices to show that
$Y$ has $A_1$ or $A_2$ singularities at the generic point of $Y_{sing}$.  
\end{rema}

\begin{theo} \label{theo:divisorcone}
Let $(X,g)$ be a polarized irreducible holomorphic symplectic manifold.
The cone of effective divisors $\NE^1(X)$ satisfies
$$ 
\left<R:\left(R,R\right) > 0, R\cdot g >0 \right>
\subset \NE^1(X)
\subset
\left<R: \left(R,R\right) \ge -2, R\cdot g >0 \right>,$$
where $R$ is taken in $\rN_1(X,\bZ)$. 
\end{theo}
Note that the cones naturally sit in
$\rN_1(X,\bR)$.  However, the Beauville-Bogomolov form allows
us to identify $\rN_1(X,\bR)=\rN^1(X,\bR)$.  
\begin{proof}
If $D \in \rN^1(X,\bZ)$ satisfies $\left(D,D\right)>0$ then
$D$ or $-D$ is big \cite[3.10]{HuyINV} and thus in the effective cone.
It remains to understand the effective divisors $D$ with  $\left(D,D\right)<0$.
(The existence of a non-trivial divisor $D$ on $X$ with $\left(D,D\right)=0$ is
conjectured to imply that $X$ is birational to an abelian fibration \cite[Conj.~3.8]{HTGAFA99}
\cite[Conj.~4.2]{Saw04}.)  

We are grateful to E. Markman for drawing our attention to the following
result, which combines \cite[Prop. 1.4]{Druel} and \cite[Thm. 4.5]{Bouck}:
\begin{quote}
Let $(X,g)$ be a polarized irreducible holomorphic symplectic manifold.
Suppose that $E$ is an irreducible divisor on $X$ with
$\left(E,E\right)<0$.  Then there exists a smooth irreducible
holomorphic symplectic variety $X'$ birational to $X$ such that
the corresponding divisor $E'\subset X'$ is contractible.
\end{quote}
Thus our analysis of extremal contractions applies to $E' \subset X'$;  we
can use the birational map to translate this information back to $X$.  
Indeed, this map gives a natural identification
$\rH^2(X,\bZ)=\rH^2(X',\bZ)$ preserving the Beauville-Bogomolov form;
by duality, there is an induced identification
$\rH_2(X,\bZ)=\rH_2(X',\bZ)$.  

The analysis of Theorem~\ref{theo:divisorial} (extended via Remark~\ref{rema:footnote})
shows that there
exists a class $R' \in \rN_1(X',\bZ)$ proportional to $E'$ such
that $\left(R',R'\right) \ge -2$.  Now $E$ and $E'$ are identified,
as are $R$ and $R'$.  In particular, $E$ is proportional
to a class $R \in \rN_1(X,\bZ)$ with $\left(R,R\right) \ge -2$.  
\end{proof}

\begin{coro}
Let $M$ be a divisor class such that $M\cdot R\ge 0$ for each
$R \in \rN_1(X,\bZ)$ with $R\cdot g>0$ and $\left(R,R\right) \ge -2$.
Then $M$ is in the closure of the moving cone.  
\end{coro}
The closure of the moving cone coincides with the intersection
of the closure of the birational K\"ahler cone with $\rN^1(X,\bR)$
\cite[Thm. 7, Cor. 19]{HTGAFA08}.  
By \cite[Prop. 4.2]{Huy03}, a divisor $M$ is in the closure of the birational
K\"ahler cone if $\left(M,D\right) \ge 0$ for all uniruled divisors $D$.
This condition follows from Theorem~\ref{theo:divisorcone}.

\begin{rema} \label{rema:Markman}
Markman \cite[Thm. 1.2]{Mark09} establishes further intersection-theoretic
properties of exceptional divisors on varieties $X$ deformation equivalent to
length $n\ge 2$ Hilbert schemes of K3 surfaces.  His approach is to classify possible
`Picard-Lefschetz' type reflections on the cohomology lattice.  If $\rho$ is a
primitive divisor class associated to such a divisor then
$$\left(\rho,\rho\right)=-2,-2(n-1).$$
In the latter case, $\rho$ is divisible by $(n-1)$ in $\rH_2(X,\bZ)$, i.e.,
$$\rho=(n-1)m R$$
for some positive integer $m$.  Thus we have
$$\left(R,R\right)=-2,-\frac{2}{m^2(n-1)}.$$ 
\end{rema}

\section{Data for Hilbert schemes of K3 surfaces} \label{sect:Hilbtable}
In this section, we tabulate data for extremal rays with prescribed intersection
properties, with references to explicit varieties exhibiting each type.
We focus on $S^{[n]}$ for $n \le 4$.
Our objective is to provide evidence for the conjectures in the introduction.

\begin{center}
\begin{tabular}{|c|c|l|}
\hline
  & & \\
$\left(R,R\right)$ & $\left(\rho,\rho\right)$  & Geometry \\
\hline 
 & &  \\
$-2$ & $-2$ &  $\bP^1$   \\
 & &   \\
\hline
\end{tabular}

\

Table H1: $X\sim S$
\end{center}

\begin{center}
\begin{tabular}{|c|c|l|l|}
\hline
 & & & \\
$\left(R,R\right)$ & $\left(\rho,\rho\right)$  & Geometry & Example \\
\hline 
 & & & \\
$-\frac{1}{2}$ & $-2$ & Diagonal $\bP^1$-bundle over
$M\sim S$ & \ref{exam:FM}, \ref{exam:diagonal} \\
 & & & \\
$-2$ & $-2$ &  $\bP^1$-bundle over  $M\sim S$ &  \ref{exam:minustwo}\\
 & &  & \\
$-\frac{5}{2}$ & $-10$ & $\bP^2$  (Lagrangian)  & \ref{exam:degreeten}, \ref{exam:Lagr}\\
 & & & \\
\hline
\end{tabular}

\

Table H2: $X\sim S^{[2]}$
\end{center}

See \cite{HTGAFA99} for additional representative examples.  
In \cite{HTGAFA08} we proved that every extremal ray with negative
self-intersection fits into one of these three classes.  
In particular, the divisors asserted to be ample
in Conjecture~\ref{conj:main1} are indeed ample.  
In \cite{HTJIMJ} we proved the Conjecture~\ref{conj:main1}
in a specific example, i.e., the variety of lines on a cubic fourfold 
admitting a hyperplane section with six double points. 

\begin{center}
\begin{tabular}{|c|c|l|l|}
\hline
 & & & \\
$\left(R,R\right)$ & $\left(\rho,\rho\right)$ & Geometry & Example\\
\hline 
 & & & \\
$-\frac{1}{4}$ & $-4$ & Diagonal $\bP^1$-bundle over 
& \ref{exam:diagonal} \\
 & & $(S\times S)_{\circ}$ &  \\
 & & & \\
$-2$ & $-2$ & $\bP^1$-bundle over $M\sim S^{[2]}_{\circ}$ & \ref{exam:Pic} ($n=3$ case),  \ref{exam:minustwo}\\
 & & & \\
$-\frac{9}{4}$ & $-36$ & $\bP^2$-bundle over   $M\sim S$ & \ref{exam:degreefour} \\
 & & & \\
$-3$ & $-12$ & $\bP^3$   (Lagrangian)  & \ref{exam:Lagr}\\
 & & & \\
\hline
\end{tabular}

\

Table H3: $X\sim S^{[3]}$
\end{center}

\begin{center}
\begin{tabular}{|c|c|l|l|}
\hline
 & & & \\
$\left(R,R\right)$ & $\left(\rho,\rho\right)$ & Geometry & Example\\
\hline
 & & & \\
$-\frac{1}{6}$ & $-6$ & Diagonal $\bP^1$-bundle over $M\sim (S\times S^{[2]})_{\circ}$ & \ref{exam:diagonal} \\
 & & &\\
$-\frac{2}{3}$ & $-6$ & $\bP^1$-bundle over $M\sim (S^{[2]} \times S')_{\circ}$, & \ref{exam:isog} \\
 & &  where $S$
and $S'$ are isogenous &  \\
 & & &\\
$-2$ & $-2$ & $\bP^1$-bundle over $S^{[3]}$  & \ref{exam:minustwo}\\
 & &  &\\
$-\frac{13}{6}$ & $-78$ & $\bP^2$-bundle over $S^{[2]}$ & \ref{exam:degreesix} \\
 & & &\\
$-\frac{8}{3}$ & $-24$ &  $\bP^3$-bundle over $S$ & \ref{exam:degreeeight}\\
 & & &\\
$-\frac{7}{2}$ & $-14$ & $\bP^4$  (Lagrangian) & \ref{exam:Lagr}\\
 & & & \\
\hline
\end{tabular}

\

Table H4: $X\sim S^{[4]}$
\end{center}

\begin{rema}
It is natural to ask whether Grassmannians (and Grassmann bundles)
also arise as we deform
extremal rays in holomorphic symplectic manifolds.  There are 
obvious examples of embeddings
$$\Gr(m,n) \hookrightarrow X.$$
Concretely, let $S\subset \bP^3$ be a quartic K3 surface
not containing a line with polarization $f$,
$\Gr(2,4)$ the Grassmannian of lines in $\bP^3$, and
$$\begin{array}{rcl}
\Gr(2,4) & \ra & S^{[4]} \\
  L & \mapsto & L \cap S.
\end{array}
$$
However, the minimal rational curves in $\Gr(2,4)$ correspond
to lines containing a point and contained in a plane.  
These have class $f-\delta$ in $\rH_2(S^{[4]},\bZ)$.  But
these are not the only  extremal rational curves on
$S^{[4]}$.  Consider general 
4-tuples of points on $S$ that are coplanar.  These form a 
$\bP^1$-bundle over the relative Jacobian of the linear series $|f|$
on $S$.  They also have class $f-\delta$.  

There are conceptual reasons why Lagrangian submanifolds swept
out by extremal rays are projective spaces and not more general 
varieties.  If $\ell\simeq \bP^1 \subset X$ sweeps out a Lagrangian
manifold $G$ then 
$$\Def(\ell,G)=\Def(\ell,X).$$
In particular, 
$\Gamma(\Omega^1_G|\ell)=0$
thus Riemann-Roch implies
$$\chi(\Omega^1_G|\ell)=\deg(\Omega^1_G|\ell) + \dim(G) + 1 \le 0.$$
However,
this forces $G\simeq \bP^{\dim(G)}$;  see \cite{CMSB}, which
also includes classification results for isolated symplectic singularities
along these lines.  
\end{rema}

\section{Examples of extremal rays for Hilbert schemes}
\label{sect:repexam}

In general, we do not have a conceptual explanation of the intersection-theoretic
properties of negative extremal rays. 
We therefore give examples for 
each {\em type} enumerated on the table.  The monodromy analysis of 
Markman \cite{Mark06,Mark} and Gritsenko-Hulek-Sankaran \cite{GHS} show 
that orbits of primitive elements $\rho \in \rH^2(X,\bZ)$ under the monodromy representation
are governed by intersection-theoretic properties, e.g.,
by the square $\left(\rho,\rho\right)$ and the ideal 
$\left(\rho,\rH^2(X,\bZ)\right)\subset \bZ$.  

Furthermore, the deformation-theoretic results of \cite{HTGAFA99,Ran} imply that
the examples produced here are stable under deformations preserving the
Hodge class of $R$, the negative extremal ray under consideration.  

\subsection*{Review of Fourier-Mukai formalism}
Let $(S,f)$ be a complex polarized K3 surface. 
The cohomology 
$$\rH^*(S,\bZ)=\rH^0(S,\bZ) \oplus \rH^2(S,\bZ) \oplus \rH^4(S,\bZ)$$
admits a Mukai pairing
$$\left<(r_1,D_1,s_1),(r_2,D_2,s_2)\right>=D_1\cdot D_2 - r_1s_2 - r_2s_1.$$
This is a unimodular integral quadratic form.  
We regard the full cohomology ring as a Hodge structure of weight two,
with the degree zero and four parts interpreted as Tate 
Hodge structures $\bZ(-1)$.  

Let $\cE$ be a coherent sheaf on $S$.  Its Mukai vector
is defined
$$v(\cE)=(r(\cE),c_1(\cE),s(\cE)p) \in \rH^*(S,\bZ)$$
where $r(\cE)$ is the rank, $p$ is the class of a point on $S$, and 
$$s(\cE)=\chi(S,\cE)-r(\cE)=c_1(\cE)^2/2-c_2(\cE)+r(\cE).$$
Given two coherent sheaves $\cE$ and $\cF$ and $S$, we have \cite[2.2]{Muk2}
$$\chi(\cE,\cF):=\Hom(\cE,\cF)-\Ext^1(\cE,\cF)+\Ext^2(\cE,\cF)=
-\left<v(\cE),v(\cF)\right>.$$

A coherent sheaf 
$\cE$ is {\em simple} if its
only endomorphisms are homotheties, i.e., $\Hom(\cE,\cE)=\bC$.  
Mukai \cite{Muk} has shown that the moduli space of such sheaves
is smooth and holomorphic symplectic of dimension
$$2-\chi(\cE,\cE)=\left<v(\cE),v(\cE)\right>.$$

\begin{theo}
Let $(S,f)$ be a general polarized K3 surface of degree $f\cdot f$.  
Let $v$ be a Mukai vector of a simple sheaf on $S$ and $M_{v}$
the corresponding moduli space of simple sheaves.  Assume that $v$ is primitive
and isotropic with respect to the Mukai pairing.  
Then each connected component $M\subset M_v$ is a K3 surface with period 
computed by the following identification
of Hodge structures
$$\rH^2(M,\bZ)= v^{\perp}/\bZ v.$$
\end{theo}
\begin{proof}
First assume that $(S,f)$ is an arbitrary polarized K3 surface of this degree.
Given a polarization $A$ on $S$, we can consider the moduli space
$M_{A,v}$ of $A$-stable sheaves with Mukai vector $v$.  This is compact
provided the numbers $r(\cE),s(\cE),$ and $A\cdot c_1(\cE)$ are
relatively prime \cite[4.1]{Muk2}.  Since stable sheaves are
automatically simple, the moduli space of simple sheaves is also
compact.  

Provided $v$ is primitive, we can always specialize  $(S,f)$ to obtain an extra
polarization $A$ such that $M_{A,v}$ is compact.  It follows that
$M_v$ is also compact, and this remains true after deforming back to a generic
$S$.

The computation of the period of $M$ is \cite[Thm. 1.4]{Muk2}.  
\end{proof}

Let $\cE \ra S \times M$ denote a universal sheaf;
the Fourier-Mukai transform of bounded derived categories is defined \cite[1.5]{Muk2}
$$\begin{array}{rcl}
\Phi_{\cE}:D^b(S) & \ra & D^b(M) \\
           \cF & \mapsto & (\pi_M)_*(\pi_S^*\cF \otimes_{\bL} \cE),
\end{array}
$$
where we take the derived push-forward, pull-back, and tensor-product.
Let $\phi_{\cE}$ denote the induced map on cohomology, characterized
by the formula \cite[5.29]{Huy}:
$$\phi_{\cE}\left(\ch(\cF)\cdot \sqrt{\td(S)}\right)=\ch(\Phi_{\cE}(\cF))\cdot \sqrt{\td(M)}$$
where $\ch$ denotes the Chern character and $\td$ the Todd class.  
Examples of concrete computations in special cases can be found
in the work of Yoshioka, e.g., \cite[Lemma 2.1]{Yosh99}.
Our computations below are based on similar techniques.

\subsection*{Extremal rays from $\bP^1$-bundles}
\begin{exam} \label{exam:diagonal}
When $X\simeq S^{[n]}$, $\delta^{\vee}$ is naturally an extremal ray
of $S^{[n]}$.  Consider the contraction to the symmetric product
$$\begin{array}{rcl}
\gamma:S^{[n]} &\ra &  S^n/\fS_n\\
   \Sigma & \mapsto &  \mathrm{supp}(\Sigma)
\end{array}
$$
assigning to each subscheme its support.  
Consider the open set $\Delta^{[n]}_{\circ}\subset \Delta^{[n]}$ 
corresponding to
subschemes whose support consists of $n-1$ distinct points.  The
contraction $\gamma$ induces a morphism forgetting the nontrivial
scheme structure:
$$\begin{array}{rcl}
\gamma:\Delta^{[n]}_{\circ} & \ra & (S^{[n-2]} \times S)_{\circ} \\
\Sigma & \mapsto & \mathrm{red}(\Sigma).
\end{array}
$$
The fibers are isomorphic to $\bP^1$ and have class $\delta^{\vee}$
by the formula
$$\cO_{S^{[n]}}(\Delta^{[n]})|\ell  \simeq \cO_{\bP^1}(-2).$$
As we've seen, $\left(\delta^{\vee},\delta^{\vee}\right)=-1/2(n-1)$.
\end{exam}

\begin{exam} \label{exam:FM}
Consider the case where $(S,f)$ is a general K3 surface of
degree $2(r^2+r)$ for some 
positive integer $r$.  Consider the moduli space
$$M_v(S), \quad v=(r,f,(r+1)p)$$
parametrizing simple sheaves $\cE$ with invariants
$$c_1(\cE)=f, r(\cE)=r, c_2(\cE)=r^2+r-1, \chi(\cE)=2r+1.$$
The resulting isogenous K3 surface $M$ has degree $2(r^2+r)$ as well.  

We apply the Fourier-Mukai
functor to $S^{[2]}$ which is isomorphic to 
$$M_w(S), \quad w=(1,0,-p).$$
We compute
$$\phi_{\cE}(w)=((r+1)+ f' + rp')-(r+f'+(r+1)p')=1-p',$$
i.e., $\phi_{\cE}(w)$ has the invariants of ideals
defining elements of $M^{[2]}$.  
This is symmetric in $S$ and $M$ and induces a birational map
$$S^{[2]} \stackrel{\sim}{\dashrightarrow} M^{[2]}.$$
The diagonal in $M^{[2]}$ yields a divisor in $S^{[2]}$
ruled in $\bP^1$s over $M$.

We can interpret this divisor in $\rH^2(S^{[2]},\bZ)$:  The
class $\delta'=2f-(2r+1)\delta$ satisfies
$$\left(\delta',\delta'\right)=-2, \left(\delta',\rH^2(S^{[2]},\bZ)\right)=2\bZ$$
so our conjecture predicts $2\delta'$ is ruled by $\bP^1$s over a K3 surface.  
\end{exam}

\begin{exam} \label{exam:isoggeneral} Here we give examples of $\bP^1$-bundles over 
$S_{8n-16}^{[n-2]}\times S_{2n-4}$ in $S_{8n-16}^{[n]}$.  
Let $(S_{8n-16},f_{8n-16})$ be a degree $(8n-16)$ K3 surface
in $\bP^{4n-7}$.  Consider the moduli space of rank-two simple
sheaves over $S_{8n-16}$ with Chern classes $c_1(\cE)=f_{8n-16}$ and
$c_2(\cE)=(2n-4)p$.  These have invariants
$$\chi(\cE)=2(n-1), \quad v=(2,f_{8n-16},2n-4).$$
The moduli space $M_v(S_{8n-16})$ is isomorphic to a K3 surface
$S_{2n-4}$ of degree $2n-4$. 

Let $\cE$ denote the universal sheaf over $S_{8n-16}\times M_v(S_{8n-16})$
and $\cI$ the universal ideal sheaf over $S_{8n-16} \times S_{8n-16}^{[n-2]}$.
Consider the Fourier-Mukai transform defined over
$$M_v(S_{8n-16}) \times S_{8n-16} \times S_{8n-16}^{[n-2]}$$
by the rule
$$\cF=\bR{\pi_{13}}_* \left( \pi_{12}^* \cE \otimes_{\bL} \pi_{23}^*\cI \right).
$$
For generic points $\cE_m \in M_v(S_{8n-16})$ and $p_1+\ldots+p_{n-2}
\in S_{8n-16}^{[n-2]}$ we have
$$\cF_{(m,p_1+\ldots+p_{n-2})}=\Gamma(S_{8n-16},\cE_m(-p_1-\ldots-p_{n-2}))$$
which is two-dimensional as
$$\chi(\cE_m \otimes \cI_{p_1+\ldots + p_{n-2}})=\chi(\cE_m)-2(n-2)=2;$$
see \cite[\S 2]{Yosh99a} (esp. Lemma 2.6) for a general discussion of jumping
behavior for dimensions of spaces of global sections over moduli
spaces of sheaves.

Thus we obtain the diagram
$$\begin{array}{rcccl}
\bP^1 & \hookrightarrow & \bP(\cF) & \stackrel{c_2}{\ra} & S_{8n-16}^{[n]} \\
      &                 &{\scriptstyle \alpha} \downarrow \quad & & \\
      &                 & M_v(S_{8n-16}) \times S_{8n-16}^{[n-2]} & & 
\end{array}
$$
which gives the desired fibered subvariety.
\end{exam}

\begin{exam} \label{exam:isog}
We specialize this analysis to $n=4$, 
in which case we have Mukai-isogenous
K3 surfaces $S_{16} \subset \bP^9$ and $S_4 \subset \bP^3$.
These have been extensively studied \cite{IR}\cite{Mukg9}.  

We claim that the fibers $\alpha$ map to
smooth rational curves $\bP^1\subset S_{8n-16}^{[n]}$
with class 
$$R=f_{16}-10\delta^{\vee},$$
which satisfies $\left(R,R\right)=-2/3$.  

Choose a generic sheaf $\cE_m$ over $S_{16}$ corresponding to a point
$m \in M_v(S_{16})\simeq S_4$.  Recall that $c_1(\cE_m)=f_{16}$,
$c_2(\cE_m)=6$, and $\chi(\cE_m)=6$.  Since $\cE_m$ is generically
globally generated by six sections, we obtain a classifying map
$$\begin{array}{rcl}
S_{16} & \dashrightarrow & \Gr(4,6) \\
s & \mapsto & \Lambda_s,
\end{array}$$
uniquely defined up to projectivities.  
Choose $s_1+s_2 \in S^{[2]}_{16}$ and corresponding codimension two
linear subspaces
$$\Lambda_1,\Lambda_2 \in \bP^5.$$
Let $\ell_{12} \subset \bP^5$ denote the line where they intersect.

Let $H \subset S_{16}$ denote the hyperplane section defined by
the Schubert class associated to $\ell_{12}$, i.e.,
$$H=\{ s\in S_{16}:\Lambda_s \cap \ell_{12} \neq \emptyset \}$$
which contains $s_1$ and $s_2$ with multiplicity $>1$.  Indeed, 
$\ell_{12}$ is contained in $\Lambda_1$ and $\Lambda_2$, not just
incident to them.  We conclude that $H$ has nodes at $s_1$ and $s_2$.
Let $\tilde{H}$ denote the normalization of $H$, which is generically
of genus $9-2=7$.

Consider the induced linear series
$$\begin{array}{rcl}
\varphi: \tilde{H} & \ra & \ell_{12} \\
    s  & \mapsto & \Lambda_s \cap \ell_{12}.
\end{array}
$$
A Schubert-class computation 
$$\sigma_2=\sigma_1^2-\sigma_{11}$$
gives 
$$\deg(\varphi)=\sigma_2\cdot S_{16}-2=16-10-2=4,$$
taking into account that $s_1,s_2$ are basepoints.  

We conclude that $R$ corresponds to a $g^1_4$ on a genus-seven
hyperplane section to $S_{16}$.  Thus $R\cdot f_{16}=16$ and
$R\cdot \delta^{\vee}=10$ by the Hurwitz formula.  
We deduce that
$$R=f_{16}-10\delta^{\vee} \in \rH_2(S^{[4]}_{16},\bZ)$$
and $\left(R,R\right)=-2/3$.  
\end{exam}

\begin{exam} \label{exam:Pic}
Let $(S_{2(n-2)},f)$ denote a polarized K3 surface.
Consider the relative Picard scheme
$$\Pic^n(S_{2(n-2)},|f|)=M_{v}(S_{2(n-2)}), \quad v=(0,f,2).$$
In other words, these are degree $n$ invertible sheaves
on hyperplane sections of $S_{2(n-2)}$, which have genus $n-1$.  
Consider the relative Hilbert scheme
$$S_{2(n-2)}^{[n]}(|f|) \subset S_{2(n-2)}^{[n]}$$
consisting of 
length $n$ subschemes $\Sigma \subset S_{2(n-2)}$
lying in hyperplane sections.  The Abel-Jacobi map
$$\alpha: S_{2(n-2)}^{[n]}(|f|) \ra \Pic^n(S_{2(n-2)},|f|)$$
is a $\bP^1$-bundle.  (Generally a degree $n$ line bundle on
a genus $n-1$ curve has two sections by Riemann-Roch.)  
Thus we obtain the diagram
$$\begin{array}{rcccl}
\bP^{1} & \hookrightarrow & S_{2(n-2)}^{[n]}(|f|) & \ra & S_{2(n-2)}^{[n]} \\
      &                 &\alpha \downarrow \quad & & \\
      &                 & M_{0,f,2}(S_{2(n-2)}) & & 
\end{array}
$$
Note that this moduli space is deformation equivalent to $T^{[n-1]}$ for $T$
a K3 surface (cf. \cite[Thm. 0.1]{Yosh99a}).  Indeed, specialize to the situation where $S_{2(n-2)}$ contains a line
$\ell\simeq \bP^1$, i.e., $f\cdot \ell =1$.  There is an isomorphism
$$\begin{array}{rcl}
\Pic^n(S_{2(n-2)},|f|) & \ra & \Pic^{n-1}(S_{2(n-2)},|f|) \\
      \cE & \mapsto & \cE \otimes \cO(-\ell)
\end{array}
$$
and the latter variety is birational to $S_{2(n-2)}^{[n-1]}$ by the cycle class map.
Any two birational holomorphic symplectic manifolds are deformation equivalent
\cite[Theorem 4.6]{HuyINV}.
\end{exam}

\begin{exam} \label{exam:minustwo}
Let $S$ be a K3 surface containing a $(-2)$-curve $E \subset S$.  Abusing notation,
we also write $E$ for the class of the divisor
\begin{equation}
\{Z: Z \cap E \neq \emptyset \} \subset S^{[n]}. \label{eq:divisor}
\end{equation}
This divisor admits an open subset isomorphic to 
$$(S\setminus E)^{[n-1]} \times \bP^1.$$
Let $R$ denote the class of a generic ruling;  we have
$R \cdot E = -2$ hence  $R$ equals $E$ (via the
inclusion $\rH^2(S^{[n]},\bZ) \subset \rH_2(S^{[n]},\bZ).$)  

The class $R$ is {\em not} extremal in the cone of curves of $S^{[n]}$;   
indeed $\Pi\simeq E^{[n]} \subset S^{[n]}$ is Lagrangian and the
class of a line $\ell$ is the extremal class (see Example~\ref{exam:Lagr}.)  Nevertheless,
we can deform $S^{[n]}$ so that $R$ remains algebraic but $[\ell]$ does not,
and thus the rulings of (\ref{eq:divisor}) deform to rational curves in nearby fibers \cite{HTGAFA99,Ran}.  
\end{exam}

\subsection*{Extremal rays from $\bP^{n-2}$-bundles}
Let $S_{4n-10} \subset \bP^{2n-4}$ denote a generic K3 surface of
degree $4n-10$ with polarization $f$.  
Let $n \ge 3$ and 
consider the moduli space of simple sheaves $\cE$ on $S_{4n-10}$
with invariants
$$c_1(\cE)=f, c_2(\cE)=n,r(\cE)=2.$$
We can compute
$$\chi(\cE)=n-1, \quad s(\cE)=\chi(\cE)-r(\cE)=n-3$$
and the Mukai vector 
$$v=(2,f,n-3), \quad 
\left<v,v \right>=(4n-10)-2\cdot 2 \cdot (n-3)=2.$$
Thus $M_v(S_{4n-10})$ has dimension four.

Let $\cE$ be the universal bundle over $S_{4n-10} \times M_v(S)$ and
$\cF=(\pi_2)_*\cE$ which is a vector bundle of rank $n-1$.  
We have the following diagram
$$\begin{array}{rcccl}
\bP^{n-2} & \hookrightarrow & \bP(\cF) & \stackrel{c_2}{\ra} & S_{4n-10}^{[n]} \\
      &                 & \downarrow & & \\
      &                 & M_v(S_{4n-10}). & & 
\end{array}
$$

\begin{exam} \label{exam:degreesix}
We specialize to the case where $n=4$, so the K3 surface 
$$S_{4n-10}=S_6 \subset \bP^4$$
has degree six.  Here we consider rank-two simple sheaves with
$c_1(\cE)=f,c_2(\cE)=4$, and $\chi(\cE)=3$.  The classifying maps 
associated to such $\cE$ are degree four rational maps
$$\mu_{\cE}:S_6 \dashrightarrow \bP^2$$
with the line class on $\bP^2$ pulling back to the polarization on $S_6$.
Let $Q$ denote the rank-two universal quotient bundle on $\bP^2$ so that
$\mu^*Q=\cE$ after extending over the indeterminacy.  

We construct the classifying maps explicitly:  Fix an ordinary secant line
$L(s_1+s_2)$ of $S_6 \subset \bP^4$
meeting $S_6$ at $s_1$ and $s_2$.  Projection
induces a rational map of degree four
$$\pi_{L(s_1+s_2)}: S_6 \dashrightarrow \bP^2,$$ 
resolved by blowing up $s_1$ and $s_2$.  The association
$$\begin{array}{rcl}
\phi:S_6^{[2]} & \dashrightarrow & M_{(2,f,1)}(S_6) \\
s_1 + s_2 & \mapsto & \pi_{L(s_1+s_2)}^*Q 
\end{array} 
$$
is birational.  Thus $S_6^{[2]}$ is deformation equivalent to the
compact moduli space associated with $M_{(2,f,1)}(S_6)$.  
Note that the map $\phi$ is not regular;  it has indeterminacy along
the three-secant lines of $S_6$, parametrized by the maximal 
isotropic subspaces of the quadric hypersurface containing $S_6$,
which is isomorphic to $\bP^3$.  
\end{exam}

\subsection*{Extremal rays from $\bP^{n-1}$-bundles}
Let $S_{4n-8} \subset \bP^{2n-3}$ be a general K3 surface
of degree $4n-8$ with polarization $f$;  assume $n\ge 3$.  
Consider the moduli space of simple sheaves on $S_{4n-8}$
with the following invariants:
$$c_1(\cE)=f, \quad c_2(\cE)=n, \quad r(\cE)=2.$$
We get the auxiliary invariants
$$\chi(\cE)=n, \quad s(\cE)=\chi(\cE)-r(\cE)=n-2$$
which determine a Mukai vector
$$v=(r(\cE),c_1(\cE),s(\cE)), \quad 
\left<v,v \right>=c_1(\cE)^2-2r(\cE)s(\cE)=0.$$
Let $M_v(S_{4n-8})$ denote the corresponding moduli
space, which is isomorphic to a K3 surface $T$.  
Let $\cE$ be the universal bundle over $S_{4n-8} \times M_v(S)$.
Let $\cF=(\pi_2)_*\cE$ which is a vector space of dimension four.  
We have the following diagram
$$\begin{array}{rcccl}
\bP^{n-1} & \hookrightarrow & \bP(\cF) & \stackrel{c_2}{\ra} & S_{4n-8}^{[n]} \\
      &                 & \downarrow & & \\
      &                 & M_v(S_{4n-8})&\simeq T & 
\end{array}
$$

\begin{exam} \label{exam:degreefour}
Consider the special case $n=3$, i.e., rank two simple sheaves $\cE$ on a quartic
surface $S_4$ with $c_1(\cE)=f,c_2(\cE)=3$.  We expect the classifying maps
of these to be rational maps
$$\mu:S_4 \dashrightarrow \Gr(1,3)\simeq \bP^2,$$
well-defined where $\cE$ is free and globally generated.  The map $\mu$ has 
degree three and $\cE$ coincides with an extension of $\mu^*Q$,
where $Q$ is the tautological quotient sheaf on $\bP^2$.  
(Note that $\mu^*Q$ is well-defined away from a finite subset of $S_4$,
and thus admits unique extension to a torsion-free sheaf
with the specified invariants.)
Thus $\mu$ is the projection from some point $p \in S_4$;
the corresponding bundle is denoted $\cE_p$.  

What are the loci where sections of $\cE_p$ vanish?  We assume
$S_4$ does not contain a line.  Sections of $Q$
vanish at points of $\bP^2$.  For each line 
$$p \in L_t \subset \bP^3, \quad t\in \bP^2,$$
we have 
$$L_t \cap S_4 = \{p,s_1,s_2,s_3 \}.$$
Scheme-theoretically, we take the residual scheme to $p$ in 
$L_t \cap S_4$, which is a well-defined length three subscheme
$\Sigma_t \subset S_4$.
Thus we get a morphism
$$\begin{array}{rcl}
\bP^2 & \ra & S^{[3]}_4\\
   t & \mapsto & \Sigma_t.
\end{array}
$$
\end{exam}

\begin{exam} \label{exam:degreeeight}
In the special case $n=4$ we recover an example of Mukai 
\cite[Example 0.9]{Muk}.  In this situation $M_v(S_8)\simeq S_2$,
a degree two K3 surface.  Our diagram takes the form
$$\begin{array}{ccc}
   \bP(\cF) & \stackrel{c_2}{\ra} & S^{[4]}_8 \\
{\scriptstyle q}   \downarrow \quad & & \\
   S_2  & & 
\end{array}
$$
where $q$ is a $\bP^3$ bundle.  
\end{exam}

\subsection*{Extremal rays from Lagrangian subspaces}
Here is a simple-sheaf construction of a
Lagrangian $\bP^n\subset S^{[n]}$ for special K3 
surfaces $S$.  
\begin{exam} \label{exam:Lagr}
Let $S$ be a K3 surface containing a nonsingular
rational curve $E \subset S$, which induces
$$\Pi \simeq E^{[n]} \subset S^{[n]}.$$
Lines in $\Pi$ correspond to pencils of binary forms of degree $n$,
e.g., families of subscheme
$$\ell=\{e_1+\ldots+e_n: e_1,\ldots,e_{n-1} \text{ fixed },
        e_n \in E \text{ varying } \}.$$
Since the discriminant of such a form has degree $2(n-1)$
in the coefficients, we find
$$\Delta^{[n]} \cdot \ell= 2(n-1).$$
Let $E$ denote the divisor in $S^{[n]}$ parametrizing
subschemes with some support on $E$.  Deforming the $e_1,\ldots,e_{n-1}$
to generic points of $S$, we find that
$\ell \cdot E = -2$ and
$$\ell=E - (n-1) \delta^{\vee}.$$
\end{exam}

\begin{exam} \label{exam:degreeten}
Let $(S_{10},f)$ be a polarized K3 surface of degree ten.  We can
express 
$$S_{10}=\Lambda \cap Q \cap \Gr(2,5)$$
where $Q$ is a quadric hypersurface and $\Lambda$ is a linear
subspace of codimension three.  The variety of lines on
$\Lambda \cap \Gr(2,5)$ is isomorphic to $\bP^2$;  each
such line is $2$-secant to $S_{10}$.  This gives us an inclusion
$$\bP^2 \hookrightarrow S_{10}^{[2]}.$$

Here is a vector bundle interpretation:  Let $\cE$ denote the rigid
vector bundle represented by the moduli space $M_{2,f,3}(S_{10})$,
which is the restriction of the rank-two tautological quotient
bundle of $\Gr(2,5)$ to $S_{10}$ via the embedding given above.  
A generic point $s\in S_{10}$ determines a codimension two subspace
$K_s \subset \Gamma(S_{10},\cE)$.  Consider the locus
$$\Pi:=\{s_1+s_2 \in S_{10}^{[2]}: \dim K_{s_1} \cap K_{s_2} \ge 2 \}.$$
This coincides with the plane constructed in the previous paragraph.  
\end{exam}

Let $S_{4n-6} \subset \bP^{2n-2}$ be a general 
degree $(4n-6)$ K3 surface
with hyperplane class $f$.  (When $n=2$ the map to $\bP^2$ is
two-to-one.)  
Consider the moduli space of simple sheaves $\cE$ on $S_{4n-6}$
with the following invariants:
$$c_1(\cE)=f, \quad c_2(\cE)=2n, \quad r(\cE)=2.$$
We get the auxiliary invariants
$$\chi(\cE)=n+1, \quad s(\cE)=\chi(\cE)-r(\cE)=n-1$$
which determine a Mukai vector
$$v=(r(\cE),c_1(\cE),s(\cE)), \quad 
\left<v,v \right>=c_1(\cE)^2-2r(\cE)s(\cE)=-2.$$
Let $M_v(S_{4n-6})=\{\text{point} \}$ denote the corresponding moduli
space and
$\cE$ the universal sheaf over $S_{4n-6} \times M_v(S)$.
Let $\cF=(\pi_2)_*\cE$ which is a vector space of dimension $n+1$.  
We have the following diagram
$$\begin{array}{rcccl}
\bP^n & \hookrightarrow & \bP(\cF) & \stackrel{c_2}{\ra} & S_{4n-6}^{[n]} \\
      &                 & \downarrow & & \\
      &                 & M_v(S_{4n-6}) & & 
\end{array}
$$
which induces our Lagrangian $\bP^n$.  

\begin{exam}
We consider the special case where $n=3$, which goes back to 
Mukai \cite{MSug}.  Here we have a degree six K3 surface $S_6$ which can 
generically be expressed as a complete intersection of a smooth quadric and a
cubic in $\bP^4$:
$$S_6= W_2 \cap W_3.$$
We construct the rigid sheaf explicitly:  Express $W_2$ as a hyperplane
section of a smooth quadric fourfolds $W'_2$ and fix an isomorphism
$W'_2 \simeq \Gr(2,4)$ to the Grassmannian.  Let $Q\ra \Gr(2,4)$ denote
the universal quotient bundle, which is globally generated by four sections.
Then $\cE$ is the restriction of $Q$ to $S_6$.

The zero-sections of $Q$ trace out maximal isotropic subspaces of $W'_2$,
which restrict to maximal isotropic subspaces of $W_2$, i.e., a family
of lines $\{L_t\}$ parametrized by $\bP^3$.  The intersections
$L_t \cap S_6$ are length three subschemes of $S_6$, at least assuming
$S_6$ does not contain a line.  Therefore, we obtain a morphism
$$\begin{array}{rcl}
\bP^3 & \ra & S_6^{[3]} \\
  t & \mapsto & L_t \cap S_6.
\end{array}
$$
\end{exam}

\section{Data for generalized Kummer varieties}

\begin{center}
\begin{tabular}{|c|c|l|l|}
\hline
 & & & \\
$\left(R,R\right)$ & $\left(\rho,\rho\right)$  & Geometry & Example \\
\hline
 & & & \\
$-\frac{1}{6}$ & $-6$ & \text{Diagonal } $\bP^1$\text{-bundle  over }
$A$ & {\rm see introduction} \\
 & & & \\
$-\frac{2}{3}$ & $-6$ &  $\bP^1$\text{-bundle over } $A$ & \ref{exam:secondKummer}\\
 & &  & \\
$-\frac{3}{2}$ & $-6$ & $\bP^2$ \text{ (Lagrangian)}  & \ref{exam:KumLag}\\
 & & & \\
\hline
\end{tabular}

\

Table K2: $X\sim K_2(A)$
\end{center}

\

\

\begin{center}
\begin{tabular}{|c|c|l|l|}
\hline
 & & & \\
$\left(R,R\right)$ & $\left(\rho,\rho\right)$  & Geometry & Example \\
\hline
 & & & \\
$-\frac{1}{2(n+1)}$ & $-2(n+1)$ & \text{Diagonal } $\bP^1$\text{-bundle  over }
$A \times A^{[n-2]}$ & {\rm see }\S\ref{sect:intro} \\
 & & & \\
$-\frac{2}{n+1}$ & $-2(n+1)$ &  $\bP^1$\text{-bundle over  }& \ref{exam:secondKummer}\\
		 &           & \text{holomorphic symplectic manifold} & \\
 & &  & \\
\vdots & \vdots & \vdots  & \vdots \\
 & &  & \\
$-\frac{n+1}{2}$ & $-2(n+1)$ & $\bP^n$ \text{ (Lagrangian)}  & \ref{exam:KumLag}\\
 & & & \\
\hline
\end{tabular}

\

Tentative Table Kn:  $X\sim K_n(A)$
\end{center}

\section{Examples of extremal rays for generalized Kummer varieties}
\subsection*{Extremal rays from $\bP^1$-bundles}
\begin{exam} \label{exam:secondKummer}
Let $C$ be a smooth projective curve of genus two and $(A,\Theta)$ its Jacobian.  For
simplicity, we assume that the N\'eron-Severi group of $A$ is generated by $\Theta$.  Consider
$K_2(A)$ with 
$$\Pic(K_2(A))=\bZ \Theta \oplus_{\perp} \bZ e,  \quad  \left(\Theta,\Theta\right)=2, \left(e,e\right)=-6.$$
Each degree three line bundle $L$ on $C$ has $h^0(C,L)=2$, and is globally generated unless
$L\simeq \omega_C(p)$ for some point $p \in C$.  We can therefore consider the subvariety
$$ F'=\{Z: Z\subset A \text{ has length three}, Z \subset \tau_a(C) \text{ for some } a\in A \} \subset A^{[3]},$$
where $\tau_a$ is translation by $a$.  Note that $F'$ is a $\bP^1$-bundle over the degree three component
of the relative Picard scheme of the collection of translates of $C$ in $A$.  Fixing a reference $j:C \hookrightarrow A$
induces a morphism of Albanese varieties $j_*:\Alb_3(C) \ra \Alb_3(A)\simeq A$,
where the last isomorphism is translation by $3\times (0)$.   Since $A=\Alb_0(C)=\Alb_0(A)$,
$j^*$ is compatible with the right actions of $A$.  On the other hand,
$$(\tau_a \circ j)_*={\tau_{3a}}_* \circ j_*.$$
Restricting to the Kummer subvariety
$$F=F' \cap K_2(A)$$
yields a $\bP^1$-bundle over 
$$\{(a,b) \in A \times A:3a+b=0  \} \simeq A.$$

Let $R$ denote the class of a generic fiber $\ell$ of this bundle.  
We know that $\ell \cdot \Theta = 2$ and $\ell \cdot e=4$, because any degree three morphism
$C\ra \bP^1$ has eight ramification points, counted with multiplicities.  It follows that
$$R=\Theta - \frac{2}{3} e$$
which satisfies
$$\left(R,R\right)=-\frac{2}{3}.$$
Let $\rho=3R$ denote the smallest positive multiple of $R$ contained in $\rH^2(K_2(A),\bZ)$;
we have
$$\left(\rho,\rho\right)=-6, \quad \left(\rho, \rH^2(K_2(A),\bZ) \right)=3\bZ.$$
\end{exam}

\begin{exam} \label{exam:generalKummer}
Let $(A,\Theta)$ be an abelian surface with $(1,\sg)$ polarization $\Theta$ and N\'eron-Severi
group generated by $\Theta$.  Consider $K_{\sg}(A)$ with
$$\Pic(K_{\sg}(A))=\bZ \Theta \oplus_{\perp} \bZ e,  \quad  \left(\Theta,\Theta\right)=2\sg-2, \left(e,e\right)=-2(\sg+1).$$
Consider degree $(\sg+1)$ line bundles $L$ on curves $C$ homologically equivalent to $\Theta$.  We
have $h^0(C,L)\ge 2$ with equality for general $C$ and $L$.  We have the locus
$$ F'=\{Z: Z\subset A \text{ has length } \sg+1, Z \subset \tau_a(C) \text{ for some } a\in A, C \in |\Theta| \}
$$
in $A^{[\sg+1]}.$
Accounting for the fibers of $C^{[\sg+1]}\ra J^{\sg+1}(C)$, translations of $C$, and linear equivalences within each translated
divisor class, we find
$$\dim(F')=1+\sg+2+\sg-2=2\sg+1.$$
In particular, $F'$ is a divisor fibered as follows:

$$\begin{array}{ccc}
\bP^1 & \ra & F' \\
      & & \downarrow \\
J_{\sg+1}(C') & \ra & B'_1 \\
       & &\downarrow  \\
 &  & A/\left<H\right> \\
\end{array}
$$
where
\begin{itemize}
\item{the first fibration reflects $g^1_{\sg+1}$'s on deformations $C'$ of $C$;}
\item{
the second corresponds to the degree $(\sg+1)$ relative Jacobian of the $C'$; and }
\item{ $H$ is the torsion subgroup preserving $\Theta$.}
\end{itemize}  
Let $F=F' \cap K_{\sg}(A)$, which is also fibered in $\bP^1$'s over a $(2\sg-2)$-dimensional 
irreducible holomorphic symplectic manifold (cf. \cite[Thm. 1.4]{Yosh99},
\cite[Thm. 0.1]{Yosh01}).  
Let $R$ be the class of a ruling which, which meets the diagonal divisor $2e$ in $4\sg$ points
and $\Theta$ in $2\sg-2$ points.  Thus
$$R=\Theta - \frac{\sg}{(\sg+1)}e$$
with 
$$\left(R,R\right)=2\sg-2-\frac{2\sg^2}{(\sg+1)}=-\frac{2}{\sg+1}.$$
Let $\rho=(\sg+1)R$, the smallest positive multiple in $\rH^2(K_{\sg}(A),\bZ)$;
then $\left(\rho,\rho\right)=-2(\sg+1)$.  Again, this is equal to $\left(e,e\right)$.  
\end{exam}

\subsection*{Extremal rays from Lagrangian subspaces}
\begin{exam}  \label{exam:KumLag}
Let $(E_1,p_1)$ and $(E_2,p_2)$ be elliptic curves.  Let $A=E_1 \times E_2$ and
$X=K_n(A)$.  Consider the projective space
$$P=\{D \times p_2: \text{ where } D \in |(n+1)p_1| \} \simeq \bP^n.$$
Let $\ell \subset P$ be a line.  Then its class 
$$R=[\ell] \in \rH_2(X,\bZ) \subset \rH^2(X,\bZ) \otimes_{\bZ} \bQ$$
can be expressed
$$R=E_1-\frac{1}{2}e.$$

In particular, 
$$\left(R,R\right)=-(n+1)/2, \quad \rho:=2R \in \rH^2(X,\bZ), \left(\rho,\rho\right)=-2(n+1).$$
This again is equal to $\left(e,e\right)$.  

Indeed, since $\Pic(K_n(A))$ is generated by the classes of $E_1,E_2,$ and $e$, so it remains to
solve for the coefficients of
$$R=c_1 E_1 + c_2 E_2 + c_0 e.$$
Interpreting $\ell$ is a pencil of degree $n+1$ on $E_1$, we have
$$\ell \cdot E_1 =0, \quad \ell \cdot E_2 = 1, \quad \ell \cdot e=(n+1).$$
\end{exam}

\bibliographystyle{plain}
\bibliography{curvetable}

\begin{thebibliography}{10}

\bibitem{beauville}
A.~Beauville.
\newblock Vari\'et\'es {K}\"ahleriennes dont la premi\`ere classe de {C}hern
  est nulle.
\newblock {\em J. Differential Geom.}, 18(4):755--782 (1984), 1983.

\bibitem{BCHM}
C.~Birkar, P.~Cascini, Chr.~D. Hacon, and J.~McKernan.
\newblock Existence of minimal models for varieties of log general type.
\newblock {\em J. Amer. Math. Soc.}, 23:405--468, 2010.

\bibitem{Bouck}
S.~Boucksom.
\newblock Divisorial {Z}ariski decompositions on compact complex manifolds.
\newblock {\em Ann. Sci. \'Ecole Norm. Sup. (4)}, 37(1):45--76, 2004.

\bibitem{CMSB}
K.~Cho, Y.~Miyaoka, and N.~I. Shepherd-Barron.
\newblock Characterizations of projective space and applications to complex
  symplectic manifolds.
\newblock In {\em Higher dimensional birational geometry ({K}yoto, 1997)},
  volume~35 of {\em Adv. Stud. Pure Math.}, pages 1--88. Math. Soc. Japan,
  Tokyo, 2002.

\bibitem{CKM}
H.~Clemens, J.~Koll{\'a}r, and Sh. Mori.
\newblock Higher-dimensional complex geometry.
\newblock {\em Ast\'erisque}, (166):144 pp. (1989), 1988.

\bibitem{Druel}
S.~Druel.
\newblock Quelques remarques sur la d\'ecomposition de {Z}ariski divisorielle
  sur les vari\'et\'es dont la premier\`ere classe de {C}hern est nulle, 2009.
\newblock arXiv:0902.1078v2.

\bibitem{GHS}
V.~Gritsenko, K.~Hulek, and G.K. Sankaran.
\newblock Moduli spaces of irreducible symplectic manifolds.
\newblock arXiv:0802.2078, to appear in {\em Compositio Mathematica}.

\bibitem{HTGAFA99}
B.~Hassett and Y.~Tschinkel.
\newblock Rational curves on holomorphic symplectic fourfolds.
\newblock {\em Geom. Funct. Anal.}, 11(6):1201--1228, 2001.

\bibitem{HTGAFA08}
B.~Hassett and Y.~Tschinkel.
\newblock Moving and ample cones of holomorphic symplectic fourfolds.
\newblock {\em Geom. Funct. Anal.}, 19(4):1065--1080, 2009.

\bibitem{HTJIMJ}
B.~Hassett and Y.~Tschinkel.
\newblock Flops on holomorphic symplectic fourfolds and determinantal cubic
  hypersurfaces.
\newblock {\em J. Inst. Math. Jussieu}, 9(1):125--153, 2010.

\bibitem{HT10}
B.~Hassett and Y.~Tschinkel.
\newblock Hodge theory and {L}agrangian planes on generalized {K}ummer
  fourfolds, 2010.
\newblock arXiv:1004.0046.

\bibitem{HuyINV}
D.~Huybrechts.
\newblock Compact hyper-{K}\"ahler manifolds: basic results.
\newblock {\em Invent. Math.}, 135(1):63--113, 1999.

\bibitem{Huy03}
D.~Huybrechts.
\newblock The {K}\"ahler cone of a compact hyperk\"ahler manifold.
\newblock {\em Math. Ann.}, 326(3):499--513, 2003.

\bibitem{Huy}
D.~Huybrechts.
\newblock {\em Fourier-{M}ukai transforms in algebraic geometry}.
\newblock Oxford Mathematical Monographs. The Clarendon Press Oxford University
  Press, Oxford, 2006.

\bibitem{IR}
A.~Iliev and K.~Ranestad.
\newblock The abelian fibration on the {H}ilbert cube of a {$K3$} surface of
  genus 9.
\newblock {\em Internat. J. Math.}, 18(1):1--26, 2007.

\bibitem{Kal}
D.~Kaledin.
\newblock Symplectic singularities from the {P}oisson point of view.
\newblock {\em J. Reine Angew. Math.}, 600:135--156, 2006.

\bibitem{Kl}
S.~L. Kleiman.
\newblock Toward a numerical theory of ampleness.
\newblock {\em Ann. of Math. (2)}, 84:293--344, 1966.

\bibitem{LP}
E.~Looijenga and C.~Peters.
\newblock Torelli theorems for {K}\"ahler {$K3$} surfaces.
\newblock {\em Compositio Math.}, 42(2):145--186, 1980/81.

\bibitem{Mark06}
E.~Markman.
\newblock Integral constraints on the monodromy group of the hyperk\"ahler
  resolution of a symmetric product of a {K}3 surface, 2006.
\newblock arXiv:math/0601304.

\bibitem{Mark}
E.~Markman.
\newblock On the monodromy of moduli spaces of sheaves on {$K3$} surfaces.
\newblock {\em J. Algebraic Geom.}, 17(1):29--99, 2008.

\bibitem{Mark09}
E.~Markman.
\newblock Prime exceptional divisors on holomorphic symplectic varieties, 2009.
\newblock arXiv:0912.4981.

\bibitem{Muk}
Sh. Mukai.
\newblock Symplectic structure of the moduli space of sheaves on an abelian or
  {$K3$} surface.
\newblock {\em Invent. Math.}, 77(1):101--116, 1984.

\bibitem{MSug}
Sh. Mukai.
\newblock Moduli of vector bundles on {$K3$} surfaces and symplectic manifolds.
\newblock {\em S\=ugaku}, 39(3):216--235, 1987.
\newblock Sugaku Expositions {\bf 1} (1988), no. 2, 139--174.

\bibitem{Muk2}
Sh. Mukai.
\newblock On the moduli space of bundles on {$K3$} surfaces. {I}.
\newblock In {\em Vector bundles on algebraic varieties ({B}ombay, 1984)},
  volume~11 of {\em Tata Inst. Fund. Res. Stud. Math.}, pages 341--413. Tata
  Inst. Fund. Res., Bombay, 1987.

\bibitem{Mukg9}
Sh. Mukai.
\newblock Curves, {$K3$} surfaces and {F}ano {$3$}-folds of genus {$\leq 10$}.
\newblock In {\em Algebraic geometry and commutative algebra, Vol.\ I}, pages
  357--377. Kinokuniya, Tokyo, 1988.

\bibitem{Nam}
Y.~Namikawa.
\newblock Deformation theory of singular symplectic {$n$}-folds.
\newblock {\em Math. Ann.}, 319(3):597--623, 2001.

\bibitem{Ran}
Z.~Ran.
\newblock Hodge theory and deformations of maps.
\newblock {\em Compositio Math.}, 97(3):309--328, 1995.

\bibitem{Saw04}
J.~Sawon.
\newblock Abelian fibred holomorphic symplectic manifolds.
\newblock {\em Turkish J. Math.}, 27(1):197--230, 2003.

\bibitem{Wier}
J.~Wierzba.
\newblock Contractions of symplectic varieties.
\newblock {\em J. Algebraic Geom.}, 12(3):507--534, 2003.

\bibitem{Yosh99}
K.~Yoshioka.
\newblock Some examples of isomorphisms induced by {F}ourier-{M}ukai functors,
  1999.
\newblock arXiv:math/9902105v1.

\bibitem{Yosh99a}
K.~Yoshioka.
\newblock Some examples of {M}ukai's reflections on {$K3$} surfaces.
\newblock {\em J. Reine Angew. Math.}, 515:97--123, 1999.

\bibitem{Yosh01}
K.~Yoshioka.
\newblock Moduli spaces of stable sheaves on abelian surfaces.
\newblock {\em Math. Ann.}, 321(4):817--884, 2001.

\end{thebibliography}

\end{document}